\documentclass[a4paper,12pt]{article}
\usepackage{amsfonts,amsmath,amssymb,amstext,amsbsy,amsopn,amsthm}
\usepackage{url}
\usepackage{graphicx}

\newcommand{\beq}{\begin{eqnarray}}
\newcommand{\eeq}{\end{eqnarray}}
\newcommand{\beqn}{\begin{eqnarray*}}
\newcommand{\eeqn}{\end{eqnarray*}}

\newcommand{\D}{\displaystyle}
\newcommand{\itg}{\int \limits}
\newcommand{\R}{{\mathbf R}}
\newcommand{\N}{{\mathbf N}}
\newcommand{\Z}{{\mathbf Z}}
\newcommand{\Ree}{\mbox{\rm Re}\:}

\newcommand{\cD}{\mathcal D}

\newcommand{\cK}{\mathcal K}        
\newcommand{\cL}{\mathcal L}        
\newcommand{\cM}{\mathcal M}        
\newcommand{\cO}{\mathcal O}        
\newcommand{\cP}{\mathcal P}

\newcommand{\cS}{\mathcal S}        
        
\newcommand{\cF}{\mathcal F}        
\newcommand{\eps}{\varepsilon}

\newcommand{\bdx}{\mathbf{x}}
\newcommand{\bdy}{\mathbf{y}}
\newcommand{\bdm}{\mathbf{m}}
\newcommand{\bdk}{\mathbf{k}}
\newcommand{\bdxi}{\boldsymbol{\xi}}
\DeclareMathOperator{\erf}{erf} 
\DeclareMathOperator{\erfc}{erfc} 
\DeclareMathOperator{\wofz}{w} 
\DeclareMathOperator{\e}{e} 
 
\def\R{{\mathbb R}}
\def\Z{{\mathbb Z}}

\newtheorem{thm}{Theorem}[section]

\newcommand{\de}{\partial}
\newcommand{\keywords}{{\bf Keywords.  }}
\newcommand{\subjclass}{{\bf  Mathematics Subject Classification (2000). }}

\begin{document}
\title{Tensor product approximations of \\high dimensional potentials}

\author{{\small Flavia Lanzara$^{\mbox{\tiny 1}}$ , Vladimir Maz'ya$^{\mbox{\tiny 2}}$,
Gunther Schmidt$^{\mbox{\tiny 3}}$}}
\date{}
\maketitle
\baselineskip=0.9\normalbaselineskip
\vspace{-7pt}
\hspace*{-5mm}
\parbox{10cm}{
\begin{flushleft}
{\footnotesize\em
\begin{itemize}
\item[$^{\mbox{\tiny\rm 1}}$]\rm Dipartimento di Matematica, Universit\`a
``La Sapienza'',\\
Piazzale Aldo Moro 2, 00185 Rome, Italy\\
\texttt{\rm lanzara\symbol{'100}mat.uniroma1.it}
\item[$^{\mbox{\tiny\rm 2}}$]Department of Mathematics, University of
Link\"oping, \\ 581 83 Link\"oping, Sweden\\
\texttt{\rm vlmaz\symbol{'100}mai.liu.se }
\item[$^{\mbox{\tiny\rm 3}}$]Weierstrass Institute for Applied Analysis and
Stochastics, \\  Mohrenstr. 39,
10117 Berlin, Germany \\
\texttt{\rm schmidt\symbol{'100}wias-berlin.de}
\end{itemize}
}
\end{flushleft}
}



\baselineskip=0.9\normalbaselineskip
\vspace{-3pt}

\begin{abstract}
The paper is devoted to the efficient
computation of high-order cubature formulas
for volume potentials obtained within the framework of approximate approximations.
We combine this approach 
with modern methods of structured tensor product approximations.
Instead of performing high-dimensional discrete convolutions
the cubature of the potentials 
can be reduced to a certain number of one-dimensional convolutions
leading to a considerable
reduction of computing resources.
We propose one-dimensional integral representions of high-order cubature formulas
for n-dimensional harmonic and Yukawa potentials, which allow
low rank tensor product approximations.
\end{abstract}

\section{Introduction}
\setcounter{equation}{0}
The construction of efficient representations of multi-variate
integral operators plays a crucial role in the numerics of higher 
dimensional problems arising in a wide range of modern applications.
Let us mention multi-dimensional integral equations and volume potentials
of elliptic and parabolic partial differential operators in $\R^n$, $n \ge 3$. 

In the present paper we study 
the combination of high-order semi-analytic cubature formulas for
volume potentials 
with modern methods of structured tensor product approximations.
The cubature formulas have been obtained in \cite{Maz,MS2} using the method
of approximate approximations, see also 
\cite {MSbook} and the reference therein.

The application of tensor product approximations to 
the approximation of volume potentials
is described for example in \cite{HK,Khor,BHK}.
The main idea is to derive accurate tensor product approximations
of the density $u$ and the kernel $g$ of the convolution integral
of the form
\begin{align*} 
u(\bdx) = \sum_{p =1}^R u_p \prod_{j=1}^{n} v_j^{(p)}(x_j)
\, , 
\end{align*} 
where $\bdx=(x_1,\ldots,x_n)$ and 
\begin{align*} 
g(\bdx) =  \sum_{q =1}^R b_q \prod_{j=1}^{n} g_j^{(q)}(x_j)
\, , 
\end{align*} 
such that the convolution integral can be approximated by one-dimensional 
convolutions
\begin{align*} 
\itg_{\R^n}g(\bdx - \bdy )u(\bdy) \, d \bdy \approx
\sum_{p,q=1}^R b_q u_p \prod_{j=1}^{n}
\itg_{\R}g_j^{(q)}(x_j-y) v_j^{(p)}(y) \, dy \, .
 \end{align*} 
Then the numerical computation of the integral 
does not require 
to perform an $n$-dimensional discrete convolution, for example,
instead one has to compute $pq$ one-dimensional discrete
convolutions, which can lead to a considerable
reduction of computing time and memory requirements,
and gives the possibility to treat real world problems.

In this note we present some variants for
the tensor product approximation of
high order cubature for harmonic, Yukawa and
heat potentials as integral operators in $\R^n$.
After a brief introduction into cubature formulas based on
approximate approximations and their error behavior,
we describe their representations as one-dimensional integrals suitable
for tensor product approximation. 
For the example of the harmonic potential in $\R^n$, $3 \le n \le 6$, we report
on numerical tests for second and fourth order formulas,
which provide estimates of the tensor rank required to
approximate the action of the potential on a basis function
with a prescribed relative error.  
Further we report on numerical tests for the Yukawa potential, which can be approximated
accurately with a very small number of tensors in a wide range of arguments.
In the final section \ref{scnheat} we study the volume heat potential in 
$\R^n$, derive approximation results and 
describe its tensor product approximation.

\section{Semi-analytic cubature formulas for  potentials}
\setcounter{equation}{0}
Here we collect some results
on high-order cubature formulas
for the volume potentials of the differential operators
$-\Delta$ and $-\Delta+a^2$, $\Ree a^2 > 0$, in $\R^n$. 
\subsection{Harmonic potentials in $\R^n$}
The harmonic potential is
the inverse of the Laplace
operator and given in $\R^n$ by the formula
    \begin{equation}\label{harpot}
   {\cL}_n u(\bdx) =  \frac{\Gamma(\frac{n}{2}-1)}
   {4 \pi^{n/2}} \itg_{\R^n} \frac{u(\bdy)}{|\bdx-\bdy|^{n-2}} 
\, d\bdy \; , \quad  n \ge 3 \; .
    \end{equation}   
The volume potential provides the unique solution of the Poisson equation
\[
     -\Delta f = u \quad \textrm{in} \; \R^n \, , \quad |f(\bdx)| \le C
      |\bdx|^{n-2} \,  \;  \mbox{as} \; |\bdx| \rightarrow \infty \, .
\] 
The theory of {\it approximate approximations} proposes semi-analytic cubature formulas
for harmonic potentials by using
quasi-interpolation of the density $u$ by
functions for which the integral operator can be taken analytically.
For example, approximate $u$ by the quasi-interpolant
    \begin{equation} \label{ch9eq2}
        u_h(\bdx) = {\cD}^{-n/2} \sum_{\bdm \in
{\Z}^n} u(h\bdm) \, \eta \Big( \frac{\bdx-h\bdm}{\sqrt{\cD}h} 
\Big) \, ,
    \end{equation}
with a suitable generating function $\eta$. 
Then the sum 
    \begin{equation}\label{conapp}
\cL_n u_h (\bdx) = \frac{h^2}{\cD^{n/2-1}} \sum_{\bdm \in
{\Z}^n} u(h\bdm) \cL_n \eta \Big( \frac{\bdx-h\bdm}{\sqrt{\cD}h} 
\Big)
    \end{equation}   
is a cubature of the harmonic potential.
If one wants to compute the harmonic potential of $u$ on the given
grid $\{h\bdk\}$ then one has to compute the discrete convolution
    \begin{equation}\label{conapp1}
\cL_n u_h (h \bdk) = \frac{h^2}{\cD^{n/2-1}}  \sum_{\bdm \in
{\Z}^n}  a_{\bdk-\bdm} \, u(h\bdm) \, , \quad
\bdk \in \Z^n \, ,
    \end{equation}   
with the coefficients
    \begin{equation}\label{conapp2}
a_{\bdk} = a_{\bdk}(\cD) = 
\cL_n \eta \Big( \frac{\bdk}{\sqrt{\cD}} \Big) 
= \frac{\Gamma(\frac{n}{2}-1)}
   {4 \pi^{n/2}} \itg_{\R^n} \frac{\eta(\bdy)}{|\bdk/\sqrt{\cD}-\bdy|^{n-2}} 
\, d\bdy \, .
    \end{equation}   
It has been shown that for sufficiently smooth and compactly supported
functions the cubature formula \eqref{conapp} provide approximations
 with the error $O(h^{2M})+ O(\e^{-\cD \pi^2} h^2)$ if 
the function $\eta$ is chosen as
    \begin{equation}\label{LagM}
\eta_{2M}(\bdx)=\pi^{-n/2} \, L_{M-1}^{(n/2)}(|\bdx|^2) \, \e^{-|\bdx|^2} \, ,
    \end{equation}   
where $L_{j}^{(\gamma)}$ are the generalized Laguerre polynomials 
    \begin{equation*} 
L_{k}^{(\gamma)}(y)=\frac{\e^{\,y} y^{-\gamma}}{k!} \, \Big(
{\frac{d}{dy}}\Big)^{k} \!
\left(\e^{\,-y} y^{k+\gamma}\right), \quad \gamma > -{\rm 1} \, .
    \end{equation*}
Additionally, there holds the analytic representation
    \begin{equation} \label{ch9eta}
{\cL}_n \eta_{2M}\Big( \frac{\bdk}{\sqrt{\cD}} \Big) = \frac{1}{\pi^{n/2}} 
\frac{\cD^{n/2-1}}{4 |\bdk|^{n-2}}
\, \gamma\Big(\frac{n}{2} -1, \frac{|\bdk|^2}{\cD}\Big)
+ \frac{\e^{-|\bdk|^2/\cD}}{\pi^{n/2}}  
   \sum_{j=0}^{M-2}  \frac{L_{j}^{(n/2-1)}(|\bdk|^2/\cD)}{4 (j+1)} \,.
     \end{equation}
Here $\gamma$ is the lower incomplete Gamma function defined by
    \begin{equation} \label{incoGamma}
\gamma(a,x)= \itg_0^x \tau^{a-1}  \e^{\, -\tau } \, d\tau \, .
    \end{equation}               
In the case $k \in \N$ 
\[
\gamma(k,x) = (k-1) ! \Big( 1 - \e^{\, -x} \sum_{j=0}^{k-1}
\frac{x^j}{j!}\Big)
 \, ,
\]
whereas for odd space dimension $n$  one can use
\[
\gamma\Big(\frac{1}{2},x\Big) = \sqrt{\pi} \erf(\sqrt{x})
\]
with the error function  $\erf$ 
and  the recurrence relation
    \begin{equation*} 
\gamma(n/2+1,x) = \frac{n}{2} \, \gamma(n/2,x) - \e^{\, -x} x^{n/2} 
    \end{equation*}               
to derive analytic expressions for  ${\cL}_n \eta_{2M}$.

The asymptotic error estimate $O(h^{2M})+ O(\e^{-\cD \pi^2} h^2)$ 
for the cubature formula \eqref{conapp}
is based on
the error estimate $O(h^{2M})+\eps$
for the quasi-interpolant \eqref{ch9eq2} if the sufficiently smooth
and decaying generating function $\eta$ is subject to the moment condition
   \begin{equation}  \label{momcon} 
    \itg _{{\R} ^n} \eta(\bdx) \, d \bdx = 1
  \; , \;                 
    \itg _{{\R} ^n} 
     \bdx^{\boldsymbol{\alpha}} \eta(\bdx) \, d \bdx = 0 , \quad
    \forall \, \boldsymbol{\alpha} \; , \; 
   1  \le |\boldsymbol{\alpha}| <  2M \, .
    \end{equation}
The saturation error $\eps$ does not converge to zero for $h \to 0$, but
because of 
\begin{align}\label{satur}
 \eps = O\Big( \max_{\bdx} \Big| \sum_{\bdk \in {\Z}^n  \backslash \{ \boldsymbol{0} \}}
\cF \eta (\sqrt{\cD} \bdk) \e^{\, 2 \pi i (\bdk,\bdx)/h}\Big| \Big)
\end{align}
can be made arbitrarily small if the parameter $\cD$ is sufficiently large.
Additionally, the harmonic potential maps
the fast oscillating saturation term \eqref{satur}
into a 
function with norm of order $O(h^2 \eps)$, which establishes the 
error estimate $O(h^{2M})+ O(\e^{-\cD \pi^2} h^2)$ for the approximation of
the harmonic potential by using quasi-interpo\-lation of the density
with the generating function \eqref{LagM}.

So in numerical computations with $\cD \ge 3$ the formulas (\ref{conapp1}-\ref{ch9eta})
behave like
high order
cubature formulas for harmonic potentials 
in $\R^n$. 
The approximation of the potential on the grid $\{h\bdk\}$
can be done by fast convolutional methods, but   
one has to store the values $\{a_\bdk\}$
for a quite large subset of $\bdk \in \Z^n$.

\subsection{Yukawa potentials in $\R^n$}
The fundamental solution
of the operator $-\Delta +a^2$, $\Ree a^2 >0$ in $\R^n$ is given as
\[
\kappa_a(\bdx)= \frac{1}{(2
  \pi)^{n/2}} \Big( \frac{|\bdx|}{a}\Big)^{1-n/2}
K_{n/2-1} (a|\bdx|)\, ,
\] 
where
$K_\nu$ is the modified Bessel function of the second kind, also
known as Macdonald function, \cite[9.6]{Abr}. 
Thus the volume potential 
\begin{align*}
\itg_{\R^n} \kappa_a(\bdx-\bdy) u(\bdy) \, d \bdy
\end{align*} 
for $n=3$ also called Yukawa potential, provides the solution of the equation
\[
     (-\Delta +a^2) f= u \quad \textrm{in} \; \R^n \, .
\] 

To derive a cubature formula for that potential
we look for a solution of
\begin{equation*}
- \Delta f +a^2 f =\e^{\, -|\bdx|^2} , \; \bdx \in \R^n \, ,
\end{equation*}
which is given as the
one-dimensional integral
\begin{align} \label{pot_advdif_exp}
f(\bdx) &=\frac{ \e^{\,-|\mathbf{x}|^2} \, a^{n/2-1}}{|2\mathbf{x}|^{n/2-1}}
\itg_0^\infty K_{n/2-1} (a r) \, 
I_{n/2-1}(2|\mathbf{x}|r)\,r\, \e^{\, -r^2}\, dr \,  ,
\end{align}
where $I_\nu$ is the modified Bessel function of the first kind, see \cite[Section 5.2]{MSbook}.

Using the known analytic expressions of $I_{n+1/2}$ and $K_{n+1/2}$ (cf. \cite{Abr})
it is possible to derive analytic formula of \eqref{pot_advdif_exp}
for odd space dimension $n$. 
In particular, if $n=3$, then 
\begin{align} \label{pot_advdif_exp3}
f(\bdx)= \frac{\sqrt{\pi}}{8} \frac{\e^{\,-|\mathbf{x}|^2}}{|\mathbf{x}|}
\Big(\wofz\Big({i}\big(\frac{a}{2}-|\mathbf{x}|\big) \Big) - 
\wofz\Big({i}\big(\frac{a}{2}+|\mathbf{x}|\big) \Big)\Big) \, ,
\end{align}
where $\wofz$ denotes the scaled complementary error function
    \begin{equation} \label{Fad} 
\wofz(z)=  \e^{\, -z^2} \, \erfc(-iz)
 =\e^{\, -z^2}\Big(1+ \frac{2 i}{\sqrt{\pi}}  \itg_{0}^{z} \e^{t^2} \,
 dt \Big)  ,
    \end{equation}
and
    \begin{equation} \label{deferfc} 
\erfc(\tau) = 1- \erf(\tau) 
    \end{equation}
is the complementary error function.

Using the representation from \cite[Theorem 3.5]{MSbook}
\[
 L_{M-1}^{(n/2)}(|\mathbf{x}|^2) 
        \, \e^{-|\mathbf{x}|^2} 
= \sum_{j=0}^{M-1}
        \frac{(-1)^j}   {j! \, 4^j} \, \Delta^j  \e^{-|\mathbf{x}|^2}  \, ,
\]
one can derive as in the case of harmonic potentials
semi-analytic cubature formulas for the Yukawa potential
with the approximation rate $O(h^{2M})+ O(\e^{-\cD \pi^2} h^2)$.

\section{Tensor product expansions of potentials acting on
Gaussians}
\setcounter{equation}{0}
To obtain a tensor product approximation 
of the  second order cubature formula 
    \begin{equation}\label{ch9harm3}  
   \frac{{\cD}h^2}{(\pi {\cD})^{n/2}}
   \sum_{\mathbf{m} \in{\Z}^n} u(h\mathbf{m}) \, \cL_n
   \big(\e^{-|\,\cdot\,|^2}\big)(\mathbf{r_m})
    \end{equation}                                     
for the harmonic potential with
    \begin{equation*}
\mathbf{r_m} = \frac{\mathbf{x}-h\bdm}{\sqrt{\cD}h} 
\quad \mbox{and} \quad {\cL}_n \big(\e^{\,-|\,\cdot\,|^2}\big)(\bdx)
= \frac{1}{4 |\bdx|^{n-2}} \gamma\Big(\frac{n}{2} -1, {|\bdx|^2}\Big)\, ,
    \end{equation*}
we use the formula obtained in \cite{MS5}
    \begin{equation*}
{\cL}_n \big(\e^{\,-|\,\cdot\,|^2}\big)(\bdx)
= \frac{1}{4}
\itg_0^\infty \frac{\e^{\, - |\bdx|^2/(1+t)}}{(1+ t)^{n/2}}
\, dt  \, ,
    \end{equation*}
which is valid for $n \ge 3$.
The simple quadrature of the integral
    \begin{equation*}
{\cL}_n \big(\e^{\,-|\,\cdot\,|^2}\big)(\bdx)
\approx \sum_{k=1}^R \omega_k \frac{\e^{\, -
    |\bdx|^2/(1+\tau_k)}}{(1+ \tau_k)^{n/2}} 
= \sum_{k=1}^R \omega_k\prod_{j=1}^{n}
\frac{\e^{\, -  x_j^2/(1+\tau_k)}}{(1+ \tau_k)^{n/2}}
    \end{equation*}
with certain quadrature  weights $\omega_k$ and nodes $\tau_k$ 
gives already a tensor product approximation.
Hence, for $\bdx=(x_1,\ldots,x_n)$ and $\bdm=(m_1,\ldots,m_n)$
\[
\cL_n \big(\e^{\,-|\,\cdot\,|^2}\big)\Big(\frac{\mathbf{x}-h\bdm}{\sqrt{\cD}h}\Big) 
\approx \frac{\cD h^2}{4 (\pi\cD)^{n/2} } \sum_{k=1}^R \omega_k \prod_{j=1}^{n}
\frac{\e^{- (x_j-h m_j)^2/(\cD h^2(1+\tau_k))}}
        {(1+\tau_k)^{n/2}}
\]
which implies that one can approximate 
\[
\cL_n u_h (h \bdk) \approx\frac{\cD h^2}{4 (\pi\cD)^{n/2} }  \sum_{\bdm \in
{\Z}^n} u(h\bdm)
\sum_{k=1}^R \omega_k \prod_{j=1}^{n}
\frac{\e^{- (k_j- m_j)^2/(\cD (1+\tau_k))}}
        {(1+\tau_k)^{n/2}} \, ,
\]
$\bdk=(k_1,\ldots,k_n)$.

To obtain a similar tensor product approximation for higher order cubature
formula we note that one obtains the same convergence order
$O(h^{2M})+ O(\e^{-\cD \pi^2} h^2)$ as in the case of generating functions \eqref{LagM}
if the density is approximated by
the sum 
\begin{align*}
\cM_M u(\bdx)= (\pi \cD)^{-n/2}
\sum_{\bdm \in \Z^n} u(h\bdm)
\prod_{j=1}^n \widetilde \eta_{2M} \Big( \frac{x_j - h m_j}{\sqrt{\cD} h}\Big) \, ,
 \end{align*}
where the generating function is the tensor product of
the one-dimensional generating functions
\begin{align} \label{Lagone}
\widetilde \eta_{2M} (x)= L_{M-1}^{(1/2)}(x^2)
\e^{\, - x^2} \, .
 \end{align}
and obviously satisfies the moment condition \eqref{momcon}. 

To get the one-dimensional integral representation
of $\cL_n \Big(\prod_{j=1}^n \widetilde \eta_{2M}\Big)$
we use the relation
\begin{align*}
\widetilde\eta_{2M} (x)&=
\sum_{k=0}^{M-1}
        \frac{(-1)^k}   {k! \, 4^k} \, 
\frac{d^{2k}}{dx^{2k}}  \e^{-x^2} =\e^{\, -x^2}\sum_{k=0}^{M-1}
        \frac{(-1)^k }   {k! \, 4^k} \,  H_{2k}(x)=\e^{\, - x^2} \sum_{k=0}^{M-1}
L_{k}^{(-1/2)}(x^2)\, ,
 \end{align*}
where $H_{k}(x)$ denotes  the Hermite polynomial
     \begin{equation*} 
 H_{k}(x)=(-1)^k \e^{\,x^2} \left(\frac{d}{dx}\right)^k
 \e^{\, -x^2}
 \, .
     \end{equation*}
Then the solution of the Poisson equation
\begin{align*}
- \Delta u(\bdx)
= \prod_{j=1}^n \widetilde \eta_{2M} (x_j)
 \end{align*}
is given by the integral
\begin{align*}
&\frac{1}{4} \, \prod_{j=1}^n \Big(\sum_{k=0}^{M-1}
         \frac{(-1)^k }   {k! \, 4^k} \, 
\frac{d^{2k}}{dx_j^{2k}}\itg_0^\infty 
\frac{\e^{- x_j^2/(1+t)}}
        {(1+t)^{1/2}} \, dt \Big) \\
&=\frac{1}{4 } \itg_0^\infty 
\prod_{j=1}^n \Big(\sum_{k=0}^{M-1}
         \frac{(-1)^k }   {k! \, 4^k} \, 
\frac{d^{2k}}{dx_j^{2k}}\e^{-x_j^2/(1+t)}\Big)
\frac{ dt } {(1+t)^{n/2}} \\
&=\frac{1}{4  } \itg_0^\infty 
\prod_{j=1}^n \e^{\, -x_j^2/(1+t)} \Big(\sum_{k=0}^{M-1}
        \frac{1}   {(1+t)^{k+1/2}} \, 
L_{k}^{(-1/2)}\Big(\frac{x^2_j}{1+t}\Big)
\Big) \, dt \, .
 \end{align*}
Again, a tensor product representation 
of this integral and consequently of the
convolution matrix 
for the high order cubature of the harmonic potential is given 
by a quadrature of the last integral.

So the problem is reduced to find efficient quadrature formulas for the
parameter dependent integrals 
\begin{equation} \label{newton2}
\begin{split}
I_1(\bdx)=&\itg_0^\infty \frac{\e^{\, - |\bdx|^2/(1+t)}}{(1+ t)^{n/2}}
\, dt =\itg_0^\infty \prod_{j=1}^n \frac{\e^{\, -x_j^2/(1+t)}}{\sqrt{1+t}} \, dt \, ,
\\
I_M(\bdx)=&\itg_0^\infty 
\prod_{j=1}^n \e^{\, -x_j^2/(1+t)} \Big(\sum_{k=0}^{M-1}
        \frac{1}   {(1+t)^{k+1/2}} \, 
L_{k}^{(-1/2)}\Big(\frac{x^2_j}{1+t}\Big)
\Big) \, dt \, .
\end{split}
\end{equation}
More precisely, one has to find a certain quadrature rule
with minimal number of summands which approximates the
integrals with prescribed error for the parameters $x_j= (k_j-m_j)/\sqrt{\cD}$ within the range
$|x_j| \le K$ and some given bound $K$.

\subsection{Quadratures}
It is well known that classical trapezoidal rule 
is exponentially converging for certain classes of integrands, for example
periodic functions and rapidly decaying functions on the real line.
For example, Poisson's summation formula yields that 
\begin{align} \label{errinf}
h \sum_{k=-\infty}^ \infty f(kh)   = 
\sum_{j=-\infty}^ \infty\hat f \big( \frac{2\pi j}{h}\big)
 \end{align}
for any sufficiently smooth function, say of the Schwarz class $\cS(\R)$.
Here  $\hat f$ is the Fourier transform
\[
\hat f(\lambda) =  \itg_{-\infty}^\infty f(x) \e^{-2\pi i x \lambda} \, dx \, .
\]
Thus, 
\[
\itg_{-\infty}^\infty f(x) \, dx - 
h \sum_{k=-\infty}^ \infty f(kh) =
\sum_{j \ne 0}\hat f \big( \frac{2\pi j}{h}\big) \, ,
\]
which indicates that by choosing special substitutions
such that the integrand transforms to a rapidly decaying function with rapidly
decaying Fourier transform, the
trapezoidal rule of step size $h$ can provide very accurate
approximations of the integral.

Here we follow a 
proposal made by J. Waldvogel \cite{Wald} to compute
accurately integrals of analytic functions.
We make the substitutions
\begin{align} \label{subwald}
t=\e^{\xi}\, , \quad \xi=a(\tau + \e^{\tau}) \quad \mbox{and} \quad 
\tau=b(u -\e^{-u})\, ,
 \end{align}
where $a,b >0$ are certain constants.
Then the integrals \eqref{newton2} are transformed to integrals
over $\R$ 
of doubly exponentially decaying integrands $f$, i.e.
$|f(u)| \le c \exp(-\alpha \exp(|u|))$ for $|u| \to \infty$ with
certain constants $c, \alpha > 0$.
It is known (cf. e.g. \cite{Khor}) that by suitable truncation of
the infinite sum in 
the trapezoidal rule of step size $h$
provides 
exponentially convergent numerical quadrature algorithms
with the error estimate $O(\e^{-c/h})$. 

After the substitution we have
\begin{align*}
&I_1(\bdx)= ab \itg_{-\infty}^\infty \exp\! \Big( \!-\frac{ |\bdx|^2}
{1+ \phi(u)}\Big)
\frac{(1+\e^{-u})(1 + \exp (b(u -\e^{-u})))\,\phi(u)}
{(1+ \phi(u))^{n/2}} \,du  \, ,
 \end{align*}
where we set
\begin{align*}
\phi(u)=\exp(ab(u -\e^{-u}) + a\exp(b(u -\e^{-u}))) \, .
 \end{align*}
Similarly
\begin{align*}
&I_M(\bdx)= ab \itg_{-\infty}^\infty\prod_{j=1}^n g_j(u) du \quad \mbox{with the functions}\\
&g_j(u)=(1+\e^{-u})(1 + \exp (b(u -\e^{-u})))\,\phi(u)
\exp\!\Big(\!-\frac{x_j^2}
{1+\phi(u) }\Big)\sum_{k=0}^{M-1} 
\frac{L_{k}^{(-1/2)}\Big(\frac{x^2_j}{1+\phi(u) }\Big) }
{(1+ \phi(u))^{k+1/2}} \, .
 \end{align*}
The integrals are approximated by the 
finite sum 
\begin{equation}\label{trapez}
\itg_{-\infty}^{\infty} f(u,\mathbf{x})\, du
\approx h\, \sum_{k=-N_0}^{N_1} f(h k,\mathbf{x}),\> |\mathbf{x}|\leq K.
\end{equation}

\subsection{Numerical Results}
 \subsubsection{Approximation to the integral $I_1(\mathbf{x})$}
\begin{figure}[h] 
 \centering
    \includegraphics[width=1.9in,height=1.2in]{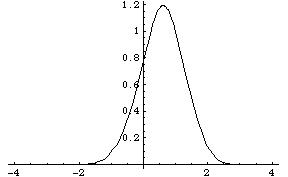} 
    \includegraphics[width=1.9in,height=1.2in]{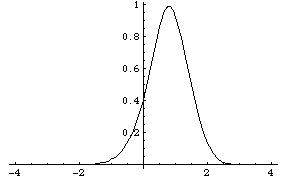}    
    \includegraphics[width=1.9in,height=1.2in]{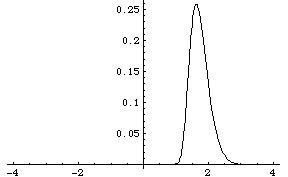}  
     \includegraphics[width=1.9in,height=1.2in]{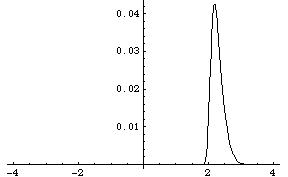} 
    \includegraphics[width=1.9in,height=1.2in]{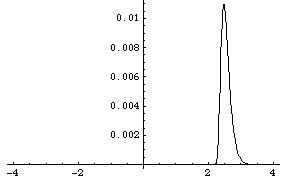}    
    \includegraphics[width=1.9in,height=1.2in]{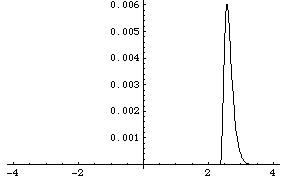}  
     \caption{The plot of the integrand function $f(u,\mathbf{x})$ ($a=b=1$) in 
$I_1(\mathbf{x})$ for $|\mathbf{x}|=0, 1, 10, 100, 500, 1000$ 
(from  the left to the right) in the interval $u\in (-4,4)$.}   
 \label{n3_ab1bis}
 \end{figure}

We assume  in  \eqref{subwald} $a=b=1$.  
Figure   \ref{n3_ab1bis}  illustrates the 
graph of  the integrand function $f(u,\mathbf{x}),\, u \in (-4,4)$, $n=3$,
 for different values of $|\mathbf{x}|\leq 10^3$.  
 A similar  behavior  holds for $n =4,5,6$.  

Table \ref{table1} presents  the maximum step $h_0$ 
and the minimum number of quadrature points required to achieve  the relative error $\epsilon$,
uniformly in  $|\mathbf{x}|\in [0,10^3]$. We have considered the space dimension  
$n=3\,, 4\,, 5\,, 6$.  

 \begin{table}[p]
\begin{center}
\begin{tabular}{ccc}\\
$n=3$ && $n=4$\\  
\begin{tabular}{|c|c|c|}
\hline
 Relative & $h_0$& Number  of\\
Error & &  quadrature points\\
\hline
$10^{-1}$ & $ 0.264 $ & $18 $\\
 \hline
$ 10^{-3}$ & $0.137 $ & $38$\\
 \hline
$ 10^{-5}$ & $0.072 $ & $82$\\
 \hline
 $10^{-7}$ & $0.055$ & $116$\\
 \hline
 $10^{-9}$ & $0.043 $ & $161$\\
 \hline
 $10^{-11} $& $0.036 $ & $205$\\
 \hline
\end{tabular}
&& 
\begin{tabular}{|c|c|c|}
\hline
 Relative & $h_0$& Number of \\
  Error & & quadrature points \\
\hline
$10^{-1}$ & $0.198$ & $20$\\
 \hline
 $10^{-3}$ & $0.095$ & $52$\\
 \hline
 $10^{-5}$ & $0.072$ & $77$\\
 \hline
 $10^{-7}$ & $0.051$ & $121$\\
 \hline
 $10^{-9} $& $0.040$ & $164$\\
 \hline
 $10^{-11}$ & $0.033$ & $206$\\
 \hline
\end{tabular}
\end{tabular}
\begin{tabular}{ccc}\\
$n=5$ && $n=6$\\  
\begin{tabular}{|c|c|c|}
\hline
 Relative  & $h_0$& Number of \\
 Error & & quadrature points \\
\hline
$10^{-1}$ & $0.181$ & $21 $\\
 \hline
 $10^{-3}$ & $0.088 $ & $59 $\\
 \hline
$10^{-5}$ & $0.065 $ & $83 $\\
 \hline
 $10^{-7}$ & $0.060$ & $96$\\
 \hline
 $10^{-9}$ & $0.037$ & $169$\\
 \hline
 $10^{-11}$ & $0.033$ & $200$\\
 \hline
\end{tabular}
&& 
\begin{tabular}{|c|c|c|}
\hline
Relative  & $h_0$& Number of \\
 Error & & quadrature points \\
\hline
$10^{-1}$ & $0.156 $ & $26 $\\
 \hline
$ 10^{-3}$ & $0.090$ & $55$\\
 \hline
 $10^{-5}$ & $0.059 $ & $90 $\\
 \hline
$ 10^{-7}$ & $0.044 $ & $130$\\
 \hline
 $10^{-9}$ & $0.035 $ & $178$\\
 \hline
$ 10^{-11}$ & $0.029 $ & $220$\\
 \hline
\end{tabular}
\end{tabular}
\caption{The approximation of $I_1(\mathbf{x})$ for $|\mathbf{x}|\leq 10^3$, with $a=b=1$ 
in \eqref{subwald}.\label{table1}}
\begin{tabular}{ccc}\\
$n=3$ && $n=4$\\  
\begin{tabular}{|c|c|c|}
\hline
 Relative & $h_0$& Number  of\\
Error & &  quadrature points\\
\hline
$10^{-1}$ & $0.0297$ & $10 $\\
 \hline
$ 10^{-3}$ & $0.0125$ & $28$\\
 \hline
$ 10^{-5}$ & $0.0077$ & $61$\\
 \hline
 $10^{-7}$ & $0.0055$ & $111$\\
 \hline
 $10^{-9}$ & $0.0042$ & $170$\\
 \hline
 $10^{-11} $& $0.0034$ & $247$\\
 \hline
\end{tabular}
&& 
\begin{tabular}{|c|c|c|}
\hline
 Relative & $h_0$& Number of \\
  Error & & quadrature points \\
\hline
$10^{-1}$ & $0.0234$ & $10 $\\
 \hline
 $10^{-3}$ & $0.0107$ & $30 $\\
 \hline
 $10^{-5}$ & $0.0070$ & $58$\\
 \hline
 $10^{-7}$ & $0.0049$ & $107  $\\
 \hline
 $10^{-9} $& $0.0037$ & $169 $\\
 \hline
 $10^{-11}$ & $0.0033 $ & $217 $\\
 \hline
\end{tabular}
\end{tabular}
\begin{tabular}{ccc}\\
$n=5$ && $n=6$\\  
\begin{tabular}{|c|c|c|}
\hline
 Relative  & $h_0$& Number of \\
 Error & & quadrature points \\
\hline
$10^{-1}$ & $0.0380$ & $7$\\
 \hline
 $10^{-3}$ & $0.0120$ & $27$\\
 \hline
$10^{-5}$ & $0.0069$ & $57$\\
 \hline
 $10^{-7}$ & $0.0046$ & $112$\\
 \hline
 $10^{-9}$ & $0.0034$ & $179$\\
 \hline
 $10^{-11}$ & $0.0031 $ & $221$\\
 \hline
\end{tabular}
&& 
\begin{tabular}{|c|c|c|}
\hline
Relative  & $h_0$& Number of \\
 Error & & quadrature points \\
\hline
$10^{-1}$ & $0.0185$ & $12$\\
 \hline
$ 10^{-3}$ & $ 0.0083$ & $ 36$\\
 \hline
 $10^{-5}$ & $0.0058$ & $70 $\\
 \hline
$ 10^{-7}$ & $0.0042$ & $117 $\\
 \hline
 $10^{-9}$ & $0.0037$ & $158 $\\
 \hline
$ 10^{-11}$ & $0.0028$ & $242$\\
 \hline
\end{tabular}
\end{tabular}
\caption{The approximation of $I_1(\mathbf{x})$ for $|\mathbf{x}|\leq 10^3$, with the choice $a=6; b=5$ 
in \eqref{subwald}.\label{table2}}
\end{center}
\end{table}%

\begin{figure}[ht] 
\centering
    \includegraphics[width=1.9in,height=1.2in]{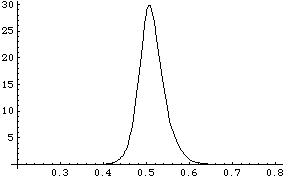} 
    \includegraphics[width=1.9in,height=1.2in]{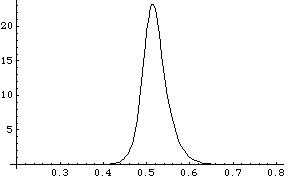}    
    \includegraphics[width=1.9in,height=1.2in]{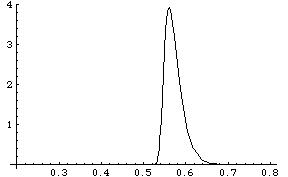}  
     \includegraphics[width=1.9in,height=1.2in]{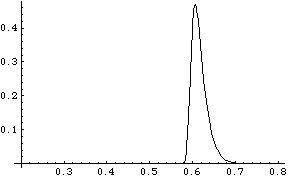} 
    \includegraphics[width=1.9in,height=1.2in]{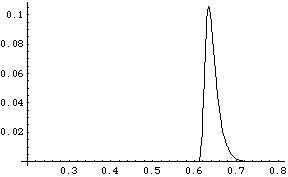}    
    \includegraphics[width=1.9in,height=1.2in]{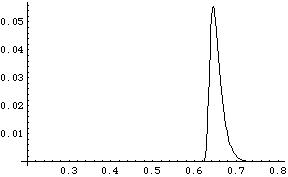}  
     \caption{The plot of the integrand function $f(u,\mathbf{x})$ ($a=6, b=5$) in $I_1(\mathbf{x})$ for 
     $|\mathbf{x}|=0, 1, 10, 100, 500, 1000$ (from  the left to the right) in the interval $u\in (0,0.85)$.}   
      \label{n3_a6b5bis}
\end{figure}

It is possible to play with different parameters  $a$ and $b$  in order
to diminish the number of summands  in the quadrature formula. Consider e.g. the case
$a=6$ and $b=5$. Figure \ref{n3_a6b5bis} shows the graph of  $f(u,|\mathbf{x}|)$, $u\in (0,0.85)$
for different values of $|\mathbf{x}|$. 
The numerical results for this quadrature are given in  Table \ref{table2}.

 \subsubsection{Approximation to the integral  $I_2(\mathbf{x})$}

  Next, we discuss the computation of the integral 
 \[
 I_2(\mathbf{x})=\itg_0^\infty \prod_{j=1}^n {\rm e}^{-\frac{x_j^2}{t+1}}
 \left( \frac{1}{\sqrt{t+1}}+\frac{1}{(t+1)^{3/2}} (\frac{1}{2}-\frac{x_j^2	}{t+1})\right)\, dt.
 \]
using the variable transformations \eqref{subwald} and the trapezoidal rule \eqref{trapez}, 
for $n=3$ and $n=4$. In the numerical results below, for the sake of simplicity,  we assumed 
$\mathbf{x}=(x,x,x)$, with $|\mathbf{x}|\leq 10^3$.

\begin{figure}[h] 
\centering
     \includegraphics[width=1.9in,height=1.2in]{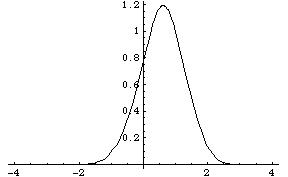} 
    \includegraphics[width=1.9in,height=1.2in]{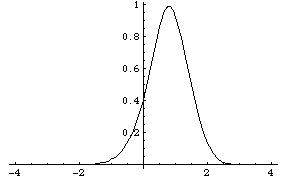}    
    \includegraphics[width=1.9in,height=1.2in]{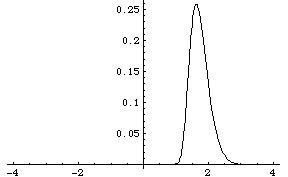}  
     \includegraphics[width=1.9in,height=1.2in]{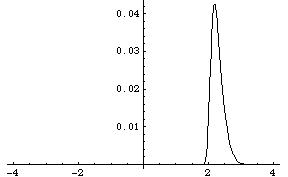} 
    \includegraphics[width=1.9in,height=1.2in]{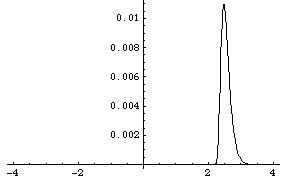}    
    \includegraphics[width=1.9in,height=1.2in]{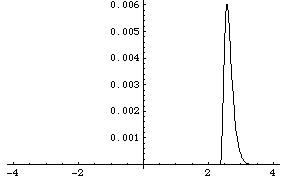} 
    \caption{The plot of the integrand function $f(u,\mathbf{x})$ in $I_2(\mathbf{x})$ 
    for $|\mathbf{x}|=0, 1, 10, 100, 500, 1000$ (from the left to the right) in the interval $u\in (-4,4)$, $a=1,b=1$.}
    \label{n3_ab1_4}
 \end{figure}

Numerical results  for this quadrature are presented in Table \ref{table3}, 
with  the parameters $a=b=1$, and in Table  \ref{table4}  in  the case $a=6,\, b=5$.

\begin{table}[ht]
\begin{center}
\begin{tabular}{ccc}\\
$n=3$ && $n=4$\\  
\begin{tabular}{|c|c|c|}
\hline
 Relative & $h_0$& Number  of\\
Error & &  quadrature points\\
\hline
$10^{-1}$ & $0.295 $ & $16 $\\
 \hline
$ 10^{-3}$ & $0.133 $ & $40$\\
 \hline
$ 10^{-5}$ & $0.072 $ & $82$\\
 \hline
 $10^{-7}$ & $0.055$ & $118$\\
 \hline
 $10^{-9}$ & $0.043$ & $163$\\
 \hline
 $10^{-11} $& $0.036$ & $204$\\
 \hline
\end{tabular}
&& 
\begin{tabular}{|c|c|c|}
\hline
 Relative & $h_0$& Number of \\
  Error & & quadrature points \\
\hline
$10^{-1}$ & $0.198$ & $23 $\\
 \hline
 $10^{-3}$ & $0.095 $ & $52 $\\
 \hline
 $10^{-5}$ & $0.072$ & $77$\\
 \hline
 $10^{-7}$ & 0.051$$ & $121$\\
 \hline
 $10^{-9} $& $0.040$ & $163 $\\
 \hline
 $10^{-11}$ & $0.033$ & $206 $\\
 \hline
\end{tabular}
\end{tabular}
\caption{The approximation of $I_2(\mathbf{x})$ for $|\mathbf{x}|\leq 10^3$, with $a=b=1$.
\label{table3}}
\begin{tabular}{ccc}\\
$n=3$ && $n=4$\\  
\begin{tabular}{|c|c|c|}
\hline
 Relative & $h_0$& Number  of\\
Error & &  quadrature points\\
\hline
$10^{-1}$ & $0.0297 $ & $10 $\\
 \hline
$ 10^{-3}$ & $0.0125 $ & $30$\\
 \hline
$ 10^{-5}$ & $0.0077$ & $63$\\
 \hline
 $10^{-7}$ & $0.0052$ & $114$\\
 \hline
 $10^{-9}$ & $0.0042$ & $175$\\
 \hline
 $10^{-11} $& $0.0034$ & $234$\\
 \hline
\end{tabular}
&& 
\begin{tabular}{|c|c|c|}
\hline
 Relative & $h_0$& Number of \\
  Error & & quadrature points \\
\hline
$10^{-1}$ & $0.0197$ & $11 $\\
 \hline
 $10^{-3}$ & $0.0107$ & $30$\\
 \hline
 $10^{-5}$ & $0.0074$ & $57$\\
 \hline
 $10^{-7}$ & $0.0046$ & $120$\\
 \hline
 $10^{-9} $& $0.0037$ & $175$\\
 \hline
 $10^{-11}$ & $0.0033$ & $222$\\
 \hline
\end{tabular}
\end{tabular}
\caption{The approximation of $I_2(\mathbf{x})$ for $|\mathbf{x}|\leq 10^3$, with $a=6; b=5$.
\label{table4}}
\end{center}
\end{table}

\section{Yukawa potential}
\setcounter{equation}{0}
To derive a tensor product approximation of the second and higher
order cubature formulas for the Yukawa potential
we use the relation
\begin{align*}
\itg_{\R^n} \kappa_a(\bdx-\bdy) \e^{\, -|\bdy|^2} \, d \bdy
= \frac{1}{4}
\itg_0^\infty \frac{\e^{-a^2 t/4}\e^{\, - |\bdx|^2/(1+t)}}{(1+ t)^{n/2}}
\, dt \, ,
\end{align*}
obtained in \cite{MS5} and is valid for all $n \ge 2$, see also \cite[Theorem 6.4]{MSbook}.
Hence, an approximate  solution of
the  equation in $\R^n$
\[
-\Delta f +a^2 f = u
\]
is given by
\begin{align*}
f_h(\bdx)=&\frac{\cD h^2}{4 (\pi\cD)^{n/2} }
 \sum_{\bdm \in  {\Z}^n}u(h \bdm)
\itg_0^\infty 
\frac{\e^{-a^2 \cD h^2 t/4}\, \e^{- |\bdx-h\bdm|^2/(\cD h^2(1+t))}}
        {(1+t)^{n/2}} \, dt \, , 
 \end{align*}
which converges with the order $O(h^2)$ to $f$.

Analogously to  the case of harmonic potentials we
consider the integral
\begin{align*}
\itg_{\R^n} &\kappa_a(\bdx- \bdy)
\prod_{j=1}^n \widetilde \eta_{2M}(y_j) \, d \bdy
= \prod_{j=1}^n \Big( \sum_{k=0}^{M-1}
         \frac{(-1)^k }   {k! \, 4^k} \, 
\frac{\partial^{2k}}{\partial x_j^{2k}} \Big)\itg_{\R^n} \kappa_a(\bdx-\bdy) \e^{- |\bdy|^2}\, d \bdy \\
&= \frac{1}{4} \, \prod_{j=1}^n \Big(\sum_{k=0}^{M-1}
         \frac{(-1)^k }   {k! \, 4^k} \, 
\frac{d^{2k}}{dx_j^{2k}}\itg_0^\infty \e^{-a^2t/4}
\frac{\e^{- x_j^2/(1+t)}}
        {(1+t)^{1/2}} \, dt \Big) \\
&=\frac{1}{4  } \itg_0^\infty 
\prod_{j=1}^n \e^{\, -x_j^2/(1+t)} \Big(\sum_{k=0}^{M-1}
        \frac{\e^{-a^2t/4}}   {(1+t)^{k+1/2}} \, 
L_{k}^{(-1/2)}\Big(\frac{x^2_j}{1+t}\Big)
\Big) \, dt \, 
 \end{align*}
which is the basis of the cubature formula of the order
$O(h^{2M})+ O(\e^{-\cD \pi^2} h^2)$.

To get doubly periodic integrands for the integrals
\begin{align*}
K_1(\bdx)=&\itg_0^\infty \frac{\e^{-a^2  t/4} \e^{\, - |\bdx|^2/(1+t)}}{(1+ t)^{n/2}}
\, dt = \itg_0^\infty \e^{-a^2  t/4} \prod_{j=1}^n \frac{\e^{\, -x_j^2/(1+t)}}{\sqrt{1+t}} \, dt\\ 
K_M(\bdx)=&\itg_0^\infty\e^{-a^2  t/4} 
\prod_{j=1}^n  \sum_{k=0}^{M-1}
        \frac{\e^{\, -x_j^2/(1+t)}}   {(1+t)^{k+1/2}} \, 
L_{k}^{(-1/2)}\Big(\frac{x^2_j}{1+t}
\Big) \, dt \, .
 \end{align*}
we make the substitutions
\[
t=\exp(b(u-\exp(-u)) \, , \; b>0 \, ,
\]
and apply the trapezoidal rule to
\begin{align*}
K_1(\bdx)=&a \itg_{-\infty}^\infty \exp\! \Big( \!-\frac{ |\bdx|^2}
{1+ \phi(u)}-\frac{c^2\phi(u)}{4}\Big)
\frac{(1+\e^{-u})\,\phi(u)}
{(1+ \phi(u))^{n/2}} \, du \, ,
 \end{align*}
\begin{align*}
K_M(\bdx)=&a\itg_{-\infty}^\infty\e^{-c^2  \phi(u)/4} 
(1+\e^{-u})\,\phi(u)  \prod_{j=1}^n \exp\! \Big( \!-\frac{ x_j^2}
{1+ \phi(u)}\Big)\\
&\hskip 3cm \times
\sum_{k=0}^{M-1}
        \frac{1}   {(1+\phi(u))^{k+1/2}} \, 
L_{k}^{(-1/2)}\Big(\frac{x^2_j}{1+\phi(u)}\Big)\, dt \, .
 \end{align*}

\subsection{ Approximation to the integral $K_1(\mathbf{x})$}

We apply the quadrature formula \eqref{trapez} to the integral $K_1(\mathbf{x})$  for $n=3$,
in the cases  $a^2=0.01$,  $a^2=0.1$ (Table \ref{table5}),  $a^2=1$ and  $a^2=4$
 (Table \ref{table6}). We assumed $b=1$.

\begin{figure}[ht] 
    \includegraphics[width=1.9in,height=1.2in]{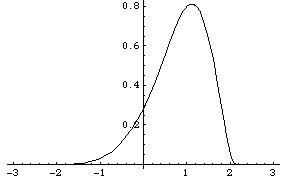} 
    \includegraphics[width=1.9in,height=1.2in]{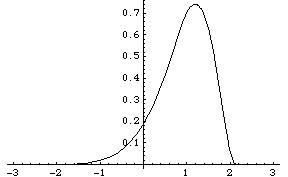}    
    \includegraphics[width=1.9in,height=1.2in]{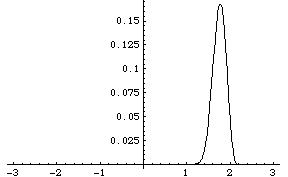}   
     \includegraphics[width=1.9in,height=1.2in]{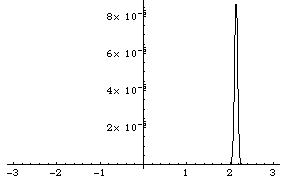} 
    \includegraphics[width=1.9in,height=1.2in]{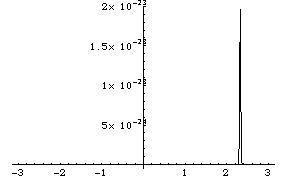}    
    \includegraphics[width=1.9in,height=1.2in]{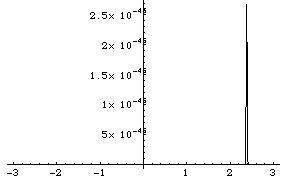}   
     \caption{The plot of the integrand function $f(u,\mathbf{x})$ in $K_1(\mathbf{x}),\> a^2=0.01, b=1$ for $|\mathbf{x}|=0, 1,10,100,500,1000$ (from  the left to the right) in the interval $u\in (-3,3)$.}   
      \label{helm1}
 \end{figure}
 
 \begin{figure}[ht] 
    \includegraphics[width=1.9in,height=1.2in]{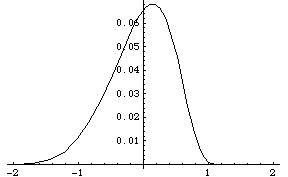} 
    \includegraphics[width=1.9in,height=1.2in]{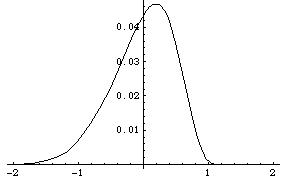}    
    \includegraphics[width=1.9in,height=1.2in]{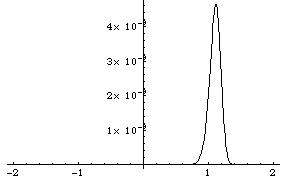}   
     \includegraphics[width=1.9in,height=1.2in]{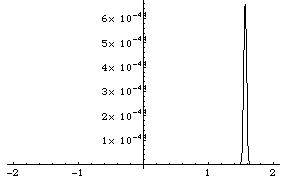}  
      \includegraphics[width=1.9in,height=1.2in]{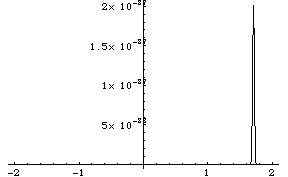}   
        \includegraphics[width=1.9in,height=1.2in]{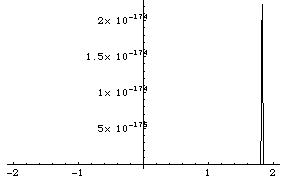}    
     \caption{The plot of the integrand function $f(u,\mathbf{x})$ in $K_1(\mathbf{x}),\> a^2=4, b=1$ 
     for $|\mathbf{x}|=0, 1, 10, 50, 100, 200 $ (from  the left to the right) in the interval $u\in (-2,2)$.}   
      \label{helm2}
 \end{figure}

 \begin{table}[h]
\begin{center}
\begin{tabular}{ccc}\\
$a^2=0.01$ && $a^2=0.1$\\  
\begin{tabular}{|c|c|c|}
\hline
 Relative & $h_0$& Number  of\\
Error & &  quadrature points\\
\hline
$10^{-1}$ & $0.99 $ & $9  $\\
 \hline
 $10^{-3}$ & $0.72 $ & $15 $\\
 \hline
 $10^{-5}$ & $0.58 $ & $20 $\\
 \hline
 $10^{-7}$ & $0.47 $ & $25 $\\
 \hline
 $10^{-9} $& $0.39 $ & $32 $\\
 \hline
 $10^{-11}$ & $0.30 $ & $43 $\\
  \hline
 $10^{-13}$ & $0.26 $ & $50 $\\
  \hline
 $10^{-15}$ & $0.25 $ & $56 $\\
\hline
\end{tabular}
&& 
\begin{tabular}{|c|c|c|}
\hline
 Relative & $h_0$& Number  of\\
Error & &  quadrature points\\
\hline
$10^{-1}$ & $0.98 $ & $7$\\
 \hline
 $10^{-3}$ & $0.68 $ & $12$\\
 \hline
 $10^{-5}$ & $0.58 $ & $17$\\
 \hline
 $10^{-7}$ & $0.42 $ & $16$\\
 \hline
 $10^{-9} $& $0.40 $ & $25$\\
 \hline
 $10^{-11}$ & $0.29 $ & $36$\\
  \hline
 $10^{-13}$ & $0.25 $ & $43$\\
  \hline
 $10^{-15}$ & $0.23 $ & $50$\\
\hline
\end{tabular}
\end{tabular}
\caption{The approximation of $K_1(\mathbf{x})$ for $|\mathbf{x}|\leq 10^3$, with $a^2=0.01$ (on the left)
and $a^2=0.1$ (on the right). \label{table5}}
\begin{tabular}{ccc}\\
$a^2=1$ && $a^2=4$\\  
\begin{tabular}{|c|c|c|}
\hline
 Relative & $h_0$& Number  of\\
Error & &  quadrature points\\
\hline
$10^{-1}$ & $0.99 $ & $6  $\\
 \hline
 $10^{-3}$ & $0.65 $ & $10 $\\
 \hline
 $10^{-5}$ & $0.48 $ & $15 $\\
 \hline
 $10^{-7}$ & $0.38 $ & $20 $\\
 \hline
 $10^{-9} $& $0.37 $ & $22 $\\
 \hline
 $10^{-11}$ & $0.29 $ & $28 $\\
  \hline
 $10^{-13}$ & $0.25 $ & $34 $\\
  \hline
 $10^{-15}$ & $0.21 $ & $42 $\\
\hline
\end{tabular}
&& 
\begin{tabular}{|c|c|c|}
\hline
 Relative & $h_0$& Number  of\\
Error & &  quadrature points\\
\hline
$10^{-1}$ & $0.92 $ & $5  $\\
 \hline
 $10^{-3}$ & $0.58 $ & $9 $\\
 \hline
 $10^{-5}$ & $0.44 $ & $13 $\\
 \hline
 $10^{-7}$ & $0.36 $ & $17 $\\
 \hline
 $10^{-9} $& $0.31 $ & $21 $\\
 \hline
 $10^{-11}$ & $0.27 $ & $25 $\\
  \hline
 $10^{-13}$ & $0.25 $ & $29 $\\
  \hline
 $10^{-15}$ & $0.16 $ & $46 $\\
\hline
\end{tabular}
\end{tabular}
\caption{The approximation of $K_1(\mathbf{x})$ for $|\mathbf{x}|\leq 10^3$, with $a^2=1$ (on the left)
and $a^2=4$ (on the right).\label{table6}}
\end{center}
\end{table}

\section{Heat potential} \label{scnheat}
\setcounter{equation}{0}
Consider the non-homogeneous (linear) heat equation
\begin{equation}\label{heatn}
 f_{t}-\nu   \Delta_{\bdx} f =u(\bdx,t),\quad \bdx\in \R^{n},\> t\geq 0
    \end{equation}
with the initial condition
\begin{equation}\label{initialn}
  f(\bdx,0)=0,\quad \bdx\in \R^{n}.
\end{equation}
It well known that 
the solution of this Cauchy problem
can be written as
\begin{equation}\label{solution}
f(\bdx,t)=\itg_0^t (\cP_{t-\lambda} u(\cdot,\lambda))(\bdx) \, d\lambda \, ,
\end{equation}
where $\cP_t$ is the Poisson integral
\begin{equation}\label{poisson}
(\cP_t	 u(\cdot,\lambda))(x)=\frac{1}{({4 \pi\nu t)^{n/2} }} \itg_{\R^n}
    {\e}^{-|\bdx-\bdxi|^2/(4\nu t)}  u(\bdxi,\lambda) \, d\bdxi \, . 
\end{equation}

An approximation of this solution $f(\bdx,t)$ can be obtained if the function $f$ is approximated
 by the quasi-interpolant on the
rectangular grid $(h\, \bdm,\tau \, j)$, with $h>0$ and $\tau>0$,
\begin{equation}\label{fht2}
    u_{h,\tau}(\bdx,t)= 
    \frac{\pi^{-(n+1)/2}}{\sqrt{\cD_{0}\cD^{n}}}\sum_{\substack{j\in\Z\\ \bdm\in\Z^{n}}}
  u(h\bdm,\tau j) \e^{-(t- \tau j)^{2}/(\cD_{0}\tau^{2})}
    \e^{-|\bdx-h\bdm|^{2}/(\cD h^{2})}.
\end{equation}
Then the sum
\begin{align}
   & f_{h,\tau}(\bdx,t)= \itg_0^t  (\cP_{t-\lambda} u_{h,\tau}(\cdot,\lambda))(\bdx)d\lambda \label{aaa}\\
&= \frac{(4 \pi\nu)^{-n/2}}{\pi^{(n+1)/2}\sqrt{\cD_{0}\cD^{n}}}
   \sum_{\substack{j\in\Z\\ \bdm\in\Z^{n}}}
    u(h\bdm,\tau j) \itg_{0}^{t}\frac{ \e^{-(\lambda- \tau j)^{2}/(\cD_{0}\tau^{2})}}
   {({t-\lambda})^{n/2}} \, d\lambda\nonumber \\
&\hskip45mm \times   \itg_{\R^{n}} {\e}^{-|\bdx-\bdxi|^{2}/(4\nu (t-\lambda))
    -|\bdxi-h\bdm|^{2}/(h^{2}\cD)} \,d\bdxi \nonumber \\
&  = \frac{h^{n}}{\pi^{(n+1)/2}\sqrt{\cD_{0}}}
    \sum_{\substack{j\in\Z\\ \bdm\in\Z^{n}}}
    u(h\bdm,\tau j) \itg_{0}^{t}\frac{ \e^{-(\lambda-(t- \tau j))^{2}/(\cD_{0} \tau^{2})}
   \e^{-|\bdx-h\bdm|^{2}/(\cD h^{2}+4\nu \lambda)}}{(\cD h^{2}+
   4\nu  \lambda )^{n/2}} \, d\lambda\nonumber
\end{align}
provides an approximation of $f(\bdx,t)$.

Since  $u_{h,\tau}(\bdx,t)$ approximates $u(\bdx,t)$ with 
\begin{align}\label{stimafn}
     \D | u(\bdx,t)-u_{h,\tau}(\bdx,t)|
   \D \leq \eps +c \big( (\tau \sqrt{\cD_{0}})^{2}+(h \sqrt{\cD})^{2}\big)\, , \quad \forall \bdx\in 
    \R^{n},\, t\in [0,T]
\end{align}
the function $f_{h,\tau}(\bdx,t)$ approximates the solution $f(\bdx,t)$ with the error 
\begin{align*}
  &  |f(\bdx,t)-f_{h,\tau}(\bdx,t)| = \D \frac{1}{({4 \pi\nu)^{n/2} }}\itg_{0}^{t} \, d\lambda
    \itg_{\R^{n}}\frac{\e^{-|\bdx-\bdxi|^{2}/(4(t-\lambda))}}{(t-\lambda)^{n/2}}
   | u(\bdxi,\lambda)-u_{h,\tau}(\bdxi,\lambda)|\, d\bdxi\\  
   & \leq T\, ||u-u_{h,\tau}||_{L^{\infty}(\R^{n} \times
  [0,T])} \leq \D \eps +c\big( (\tau \sqrt{\cD_{0}})^{2}+(h \sqrt{\cD})^{2}\big)\, , \quad \forall \bdx\in 
    \R^{n},\, t\in [0,T].
\end{align*}

The integral
\begin{equation} \label{heatint}
K_{j,\bdm}(\bdx,t)=\itg_{0}^{t}\frac{ \e^{-(\lambda-(t- \tau j))^{2}/(\cD_{0}\tau^{2})}
   \e^{-|\bdx-h\bdm|^{2}/(\cD h^{2} +4\nu \lambda)}}{(\cD h^{2} +   4\nu  \lambda )^{n/2}} \, d\lambda \\
\end{equation}
cannot be taken analytically, but it allows obviously an approximate tensor
 product approximation.
Making the substitution
\[
\lambda= \frac{t}{1+\e^{-\xi}}
\]
we derive the integral over $\R$
\begin{align*}
K_{j,\bdm}(\bdx,t)&=\frac{t}{4 }
\itg_{-\infty}^{\infty}\frac{ \e^{-(\tau j - t/(1+\e^{\xi}))^{2}/(\cD_{0}\tau^{2})}
   \e^{-|\bdx-h\bdm|^{2}/(\cD h^{2} +4\nu t/(1+\e^{-\xi}))}}
{(\cD h^{2} +   4\nu  t/(1+\e^{-\xi}) )^{n/2}\cosh^2 (\xi/2)} \, d \xi 
\end{align*}
with exponentially decaying integrand.
Performing the last 2 substitutions in \eqref{subwald}
we again transform the integrand to a doubly exponentially decaying function. 

Approximations which converge with higher order to the solution of \eqref{heatn}
can be obtained using quasi-interpolation of the right-hand side $u$ by
\begin{align} \label{fht}
\tilde u_{h,\tau}(\bdx,t) =   \frac{\pi^{-(n+1)/2}}{\sqrt{\cD_{0}\cD^{n}}}\sum_{\substack{j\in\Z\\ \bdm\in\Z^{n}}}
  u(h\bdm,\tau j)
\widetilde \eta_{2S}\Big( \frac{t - \tau j}{\sqrt{\cD_0}  \tau}\Big) 
\prod_{i=1}^n \widetilde \eta_{2M} \Big( \frac{x_i - h m_i}{\sqrt{\cD} h}\Big) 
 \end{align}
where $\widetilde \eta_{2M}$ are defined by \eqref{Lagone}.
Since for all $j \le 2S-1,\alpha_i \le 2M-1$ we have
\[
\itg_{\R^{n+1}} t^j \, \widetilde \eta_{2S}(t)
\prod_{i=1}^n x_i^{\alpha_i} \, \widetilde \eta_{2M}(x_i) \, dt \, d\bdx 
= \left\{ \begin{array}{cc} 
\pi^{(n+1)/2} \, ,& j=\alpha_1= \ldots \alpha_n=0  \, ,\\
0 \, ,& \mbox{otherwise} \, ,
\end{array} \right.
\] 
the following result can be derived in a standard way.  
\begin{thm}\label{f-fht}
    Given  $\eps>0$ there exist $\cD>0$ and $\cD_{0}>0$ such 
    that for any 
    $u\in W_{\infty}^{L}(\R^{n}\times \R)$, with $L=\max(2M,2S)$, the quasi
    interpolant \eqref{fht} satisfies the estimate
\begin{align}
|u& \,(\bdx,t) -u_{h,\tau}(\bdx,t)| \leq c_{1}(\cD h^2)^{M}+c_{2}(\cD_0\tau^2)^{S} \label{stimaf}\\
&+\eps \left( \sum_{|\alpha|=0}^{2M-1}
\frac{(h\sqrt{\cD})^{|\alpha|}}{\alpha!}||\de^{\alpha}_{x}u||_{L^{\infty}(\R^{n}\times[0,T])}+
\sum_{\beta=0}^{2S-1}\frac{(\tau
\sqrt{\cD_{0}})^{|\beta|}}{\beta!}||\de^{\beta}_{t}u||_{L^{\infty}(\R^{n}\times[0,T])}\right)\,  \nonumber
\end{align}
where the constants $c_{1}$ and $c_{2}$ do not depend on $h$, $\tau$, $\cD$,
$\cD_{0}$ and $f$.
\end{thm}

To obtain the cubature we use that
 \begin{align*}
\widetilde \eta_{2S}\Big( \frac{t}{\sqrt{\cD_0} \tau}\Big)  = \sum_{k=0}^{S}
         \frac{(-1)^k(\cD_0 \tau^2)^k }   {k! \, 4^k} \, \frac{\partial^{2k}}{\partial t^{2k}} 
 \e^{-t^{2}/(\tau^{2}\cD_{0})} \\
\widetilde \eta_{2M} \Big( \frac{x_i}{\sqrt{\cD} h}\Big)= \sum_{k=0}^{M}
         \frac{(-1)^k(\cD h^2)^k }   {k! \, 4^k} \,
\frac{\partial^{2k}}{\partial x_i^{2k}}
     \e^{-|\bdx|^{2}/(h^{2}\cD)}
\end{align*}
Hence the heat potential of the quasi-interpolant \eqref{fht}
\begin{equation} \label{heatintM}
\begin{split}
 \tilde f_{h,\tau}(\bdx,t) &= 
\itg_0^t  (\cP_{t-\lambda} \tilde u_{h,\tau}(\cdot,\lambda))(\bdx)d\lambda \\
&  = \frac{h^{n}}{\pi^{(n+1)/2}\sqrt{\cD_{0}}}
    \sum_{\substack{j\in\Z\\ \bdm\in\Z^{n}}}
    u(h\bdm,\tau j) \sum_{k=0}^{M} \cK^{S,M}_{j,\bdm}(\bdx,t) \, ,
\end{split}
\end{equation}
where we use the notation
 \begin{align*}
\cK^{S,M}_{j,\bdm}(\bdx,t)=\sum_{k=0}^{S}
         \frac{(-1)^k(\cD_0 \tau^2)^k }   {k! \, 4^k} \, \frac{\partial^{2k}}{\partial t^{2k}} 
\prod_{i=1}^{n} \sum_{k=0}^{M}
         \frac{(-1)^k(\cD h^2)^k }   {k! \, 4^k} \,
\frac{\partial^{2k}}{\partial x_i^{2k}}
K_{j,\bdm}(\bdx,t)
\end{align*}
We have
 \begin{align*}
\frac{\partial^{2k}}{\partial x_i^{2k}} \e^{-(x_i-h m_i)^{2}/(\cD h^{2} +4\nu \lambda)}
=(-1)^k k! \, 4^k\frac{\e^{-(x_i-h m_i)^{2}/(\cD h^{2} +4\nu \lambda)}}{(\cD h^{2} +4\nu \lambda)^k} 
L_{k}^{(-1/2)}\Big(\frac{(x_i-m_i)^2}{\cD h^{2} +4\nu \lambda}\Big) \, ,
\end{align*}
which by using \eqref{heatint} leads to the representation
 \begin{align*}
\cK^{S,M}_{j,\bdm}(\bdx,t) &=
\sum_{k=0}^{S}
         \frac{(-1)^k(\cD_0 \tau^2)^k }   {k! \, 4^k} \, \frac{\partial^{2k}}{\partial t^{2k}} 
\itg_{0}^{t}{ \e^{-(\lambda-(t- \tau j))^{2}/(\cD_{0}\tau^{2})}} \prod_{i=1}^{n} g_M(\lambda,x_i-m_i) d\lambda \, ,
\end{align*}
admitting again a  tensor product approximation. 
Here we denote by $g_M$ the function
 \begin{align*}
g_M(\lambda,x)= \sum_{k=0}^{M}
\frac{(\cD h^2)^k}{(\cD h^{2} +4\nu \lambda)^{k+1/2}}L_{k}^{(-1/2)}\Big(\frac{x^2}{\cD h^{2} +4\nu \lambda}\Big) 
   \e^{-x^{2}/(\cD h^{2} +4\nu \lambda)}\, .
\end{align*}

From Theorem \ref {f-fht} it is easy to deduce that 
$f_{h,\tau}$ approximates the solution $f$ with the order
$\cO((\sqrt{\cD_{0}}\tau)^{2S}+(\sqrt{\cD}
h)^{2M})$ plus the saturation error.

\begin{thm}
For any $\eps>0$ there exist $\cD>0$ and $\cD_0>0$ such that, for all 
$u \in W_{\infty}^{L}(\R^{n}\times \R)$, with $L=\max(2M,2S)$, the quasi-interpolant
\eqref{heatintM} approximates the solution of the Cauchy problem for the heat equation
\eqref{heatn}-\eqref{initialn} with the error estimate
\begin{align*}
&|f(\bdx,t)-f_{h,\tau}(\bdx,t)|\leq c_{1,T}(\cD h^2)^{M}+c_{2,T}(\cD_0\tau^2)^{S}\\&
+ \eps \left( \sum_{|\alpha|=0}^{2M-1}
\frac{(h\sqrt{\cD})^{|\alpha|}}{\alpha!}||\de^{\alpha}_{x}u||_{L^{\infty}(\R^{n}\times[0,T])}+
\sum_{\beta=0}^{2S-1}\frac{(\tau
\sqrt{\cD_{0}})^{\beta}}{\beta!}||\de^{\beta}_{t}u||_{L^{\infty}(\R^{n}\times[0,T])}\right).
\end{align*}
The constants $c_{i,T}$, $i=1,2$, depend only on $M$ and $S$.
\end{thm}
{\bf Proof.}
Since
\[
\frac{1}{({4 \pi\nu({t-\lambda}))^{n/2} }}
    \itg_{\R^{n}}
  {\e}^{-|\bdx-\bdxi|^{2}/(4\nu (t-\lambda))} d\bdxi=1
\]
we obtain that
\begin{align*}
|f(\bdx,t)&\, -f_{h,\tau}(\bdx,t)| \\
&= \frac{1}{({4 \pi\nu)^{n/2} }}\int_{0}^{t}
    \frac{d\lambda}{({t-\lambda})^{n/2}}
    \itg_{\R^{n}}
  {\e}^{-|\bdx-\bdxi|^{2}/(4\nu (t-\lambda))} |f(\bdxi,\lambda)- f_{h,\tau}(\bdxi,\lambda)| \, d\bdxi\\
  &\leq  T\, ||u-u_{h,\tau}||_{L^{\infty}(\R^{n}\times[0,T])}\, , \quad \forall \bdx\in 
    \R^{n},\, t\in [0,T].
\end{align*}
From \eqref{stimaf} the proof is complete.

\qed

\subsection{Numerical example}
We have tested the approximation formula \eqref{fht2} for solving 
the Cauchy problem
\begin{equation} \label{heatexa}
   f_{t}-  f_{xx}=x^{2}+t^{2},\qquad  f(x,0)=0, \quad x\in \R,\> t\geq 0,
    \end{equation}
having the solution
$\displaystyle{f(x,t)=t^{2}+\frac{t^{3}}{3}+tx^{2}}$.


In the  table \ref{tableheat} the
difference
$f_{h,\tau}(x,t)-f(x,t)$ for different values of $h$ and $\tau$,
$\cD_{0}=\cD=2$, at the
time $t=0.01$ and the point $x=0.01$ is given. The numerical results in the
table confirm that the error is $\cO(h^{2}+\tau^{2})$.
\bigskip

\begin{table} [h]
\begin{tabular}{|c|c|c|c|c|c|c|} 
    \hline
    $\rule[-7pt]{0mm}{20pt} \tau^{-1} \setminus h^{-1}$ & $4$ & $8$ & $16$ & $32$ & $64$ & $128$  \\
    \hline 
    $4$ & $1.25 \,10^{-3}$ &$7.81\, 10^{-4}$  &$ 6.64\, 10^{-4}$ &
    $6.35\, 10^{-4}$ & $6.27\, 10^{-4}$ & $6.25\, 10^{-4}$  \\
    \hline
    $8$ & $7.81\, 10^{-4}$ & $3.12 \,10^{-4}$ & $1.95\, 10^{-4}$ & $1.66\, 10^{-4}$
    &$1.58\, 10^{-4}$ &  $1.56\, 10^{-4}$  \\
    \hline
    $16$& $6.64\, 10^{-4}$ & $1.95 \,10^{-4}$ & $7.81\, 10^{-5}$ & $4.88\, 10^{-5}$
 & $4.15\, 10^{-5}$ & $3.96\, 10^{-5}$  \\
    \hline
    $32$ & $6.34\, 10^{-4}$ & $1.66 \, 10^{-4}$ & $4.88\,10^{-5}$ & $1.95\, 10^{-5}$ & 
    $1.22\, 10^{-5}$
    & $1.03 \, 10^{-5}$ \\
    \hline
    $64$ & $6.27\, 10^{-4}$ & $1.58\, 10^{-4}$ & $4.15\,10^{-5}$ & $1.22\,10^{-5}$ & $4.88\, 10^{-6}$ & 
    $3.05\, 10^{-6}$  \\
    \hline
    $128$ & $6.25\, 10^{-4}$ & $1.56\, 10^{-4}$ &$3.96\,10^{-5}$ & $1.03\,
    10^{-5}$ & $3.05 \, 10^{-6}$& $1.22\, 10^{-6}$  \\
    \hline
\end{tabular}
\caption{Error table for solving \eqref{heatexa} with \eqref{fht2}} \label{tableheat}
\end{table}%

\subsubsection*{Acknowledments}
The authors would like to thank B. Khoromskij
for valuable discussions concerning fast computations of 
high-dimensional problems.


\end{document}